\documentclass[11pt]{article}
\usepackage[english]{babel}
\usepackage{amssymb}
\newtheorem{lemma}{Lemma}[section]

\newtheorem{theorem}[lemma]{Theorem}
\newtheorem{proposition}[lemma]{Proposition}
\newtheorem{corollary}[lemma]{Corollary}
\newtheorem{example}[lemma]{Example}

\newtheorem{note}[lemma]{Remark}

\def\endproof{\hfill$\Box$}

\def\endproof{\hfill$\Box$}

\title{Variations of rigidity\footnote{The work
was carried out in the framework of the State Contract of the
Sobolev Institute of Mathematics, Project No.~FWNF-2022-0012.
}}
\author{S.V. Sudoplatov}

\date{}

\begin{document}

\maketitle
\begin{abstract}
We study possibilities for semantic and syntactic rigidity, i.e.,
the rigidity with respect to automorphism group and with respect
to definable closure. Variations of rigidity and their degrees are
studied in general case, for special languages and for some
natural operations with structures.

\end{abstract}


\section{Introduction}

We continue to study variations of algebraic closures
\cite{Sud_alg_23} considering and describing semantic and
syntactic possibilities for definable closures.

In Section 2, we introduce variations and degrees for semantic and
syntactic rigidity of structures, describe properties,
possibilities, and dynamics for these characteristics, in general
and for theories of unary predicates. In Section 3, indexes of
rigidity are introduced and their possibilities are described. In
Sections  4 and 5, possibilities for degrees of rigidity and for
indexes of rigidity are described for disjoint unions of
structures and for compositions of structures are studied.

We use the standard model-theoretic terminology \cite{Sh, Hodges,
PiG, Marker, TZ}, notions and notations in~\cite{Sud_alg_23}.

\section{Variations of rigidity and their characteristics}

{\bf Definition.} For a set $A$ in a structure $\mathcal{M}$,
$\mathcal{M}$ is called {\em semantically $A$-rigid} or {\em
automorphically $A$-rigid} if any $A$-automorphism $f\in{\rm
Aut}(\mathcal{M})$ is identical. The structure $\mathcal{M}$ is
called {\em syntactically $A$-rigid} if $M={\rm dcl}(A)$.

A structure $\mathcal{M}$ is called {\em $\forall$-semantically /
$\forall$-syntactically $n$-rigid} (respectively, {\em
$\exists$-se\-man\-ti\-cal\-ly / $\exists$-syntactically
$n$-rigid}), for $n\in\omega$, if $\mathcal{M}$ is semantically /
syntactically $A$-rigid for any (some) $A\subseteq M$ with
$|A|=n$.

\medskip
Clearly, as above, syntactical $A$-rigidity and $n$-rigidity imply
semantical ones, and vice versa for finite structures, but not
vice versa for some infinite ones. Besides, if $\mathcal{M}$ is
$Q$-semantically / $Q$-syntactically $n$-rigid, where
$Q\in\{\forall,\exists\}$, then $\mathcal{M}$ is $Q$-semantically
/ $Q$-syntactically $m$-rigid for any $m\geq n$.

\medskip
The least $n$ such that $\mathcal{M}$ is $Q$-semantically /
$Q$-syntactically $n$-rigid, where $Q\in\{\forall,\exists\}$, is
called the {\em $Q$-semantical / $Q$-syntactical degree of
rigidity}, it is denoted by ${\rm deg}^{Q\mbox{\rm -}{\rm
sem}}_{\rm rig}(\mathcal{M})$ and ${\rm deg}^{Q\mbox{\rm -}{\rm
synt}}_{\rm rig}(\mathcal{M})$, respectively. Here if a set $A$
produces the value of $Q$-semantical / $Q$-syntactical degree then
we say that $A$ {\em witnesses} that degree. If such $n$ does not
exists we put ${\rm deg}^{Q\mbox{\rm -}{\rm sem}}_{\rm
rig}(\mathcal{M})=\infty$ and ${\rm deg}^{Q\mbox{\rm -}{\rm
synt}}_{\rm rig}(\mathcal{M})=\infty$, respectively.

\medskip
Notice that all these characteristics have the upper bound $|M|-1$
if the structure $\mathcal{M}$ is finite. Moreover, if
$M\setminus{\rm dcl}(\emptyset)$ is finite then the cardinality
$|M\setminus{\rm dcl}(\emptyset)|-1$ is the upper bound for both
${\rm deg}^{\exists\mbox{\rm -}{\rm sem}}_{\rm rig}(\mathcal{M})$
and ${\rm deg}^{\exists\mbox{\rm -}{\rm synt}}_{\rm
rig}(\mathcal{M})$.

We have the following obvious characterizations for finite values
of degrees:

\begin{proposition}\label{prop_rig_sem_synt}
$1.$ ${\rm deg}^{\forall\mbox{\rm -}{\rm sem}}_{\rm
rig}(\mathcal{M})=0$ iff ${\rm deg}^{\exists\mbox{\rm -}{\rm
sem}}_{\rm rig}(\mathcal{M})=0$, and iff the structure
$\mathcal{M}$ is semantically rigid.

$2.$  ${\rm deg}^{\forall\mbox{\rm -}{\rm synt}}_{\rm
rig}(\mathcal{M})=0$ iff ${\rm deg}^{\exists\mbox{\rm -}{\rm
synt}}_{\rm rig}(\mathcal{M})=0$, and iff the structure
$\mathcal{M}$ is syntactically rigid.

$3.$ ${\rm deg}^{\forall\mbox{\rm -}{\rm sem}}_{\rm
rig}(\mathcal{M})=n\in\omega$ iff for any set $A\subseteq M$ with
$|A|\geq n$ there is minimal $B\subseteq A$, under inclusion, such
that $|B|=n$ and any automorphism $f\in{\rm Aut}(\mathcal{M})$
fixing $B$ pointwise fixes all elements in $\mathcal{M}$, too, and
there are no sets of cardinalities $n'<n$ with that property. Here
$B\subseteq A$ can be taken arbitrary with $|B|=n$.

$4.$ ${\rm deg}^{\exists\mbox{\rm -}{\rm sem}}_{\rm
rig}(\mathcal{M})=n\in\omega$ iff for some set $A\subseteq M$ with
$|A|\geq n$ there is minimal $B\subseteq A$, under inclusion, such
that $|B|=n$ and any automorphism $f\in{\rm Aut}(\mathcal{M})$
fixing $B$ pointwise fixes all elements in $\mathcal{M}$, too, and
there are no sets of cardinalities $n'<n$ with that property.

$5.$ ${\rm deg}^{\forall\mbox{\rm -}{\rm synt}}_{\rm
rig}(\mathcal{M})=n\in\omega$ iff for any set $A\subseteq M$ with
$|A|\geq n$ there is minimal $B\subseteq A$, under inclusion, such
that $|B|=n$ and $M={\rm dcl}(B)$, and there are no sets of
cardinalities $n'<n$ with that property. Here $B\subseteq A$ can
be taken arbitrary with $|B|=n$.

$6.$ ${\rm deg}^{\exists\mbox{\rm -}{\rm synt}}_{\rm
rig}(\mathcal{M})=n\in\omega$ iff for some set $A\subseteq M$ with
$|A|\geq n$ there is minimal $B\subseteq A$, under inclusion, such
that $|B|=n$ and $M={\rm dcl}(B)$, and there are no sets of
cardinalities $n'<n$ with that property.
\end{proposition}

By the definition, we have the following {\em monotonicity
property}: if $\mathcal{M}$ is semantically / syntactically
$A$-rigid and $A\subseteq A'\subseteq M$ then $\mathcal{M}$ is
semantically / syntactically $A'$-rigid.

Using the definition and the monotonicity property, for any
structure $\mathcal{M}$ the following inequalities hold:
\begin{equation}\label{eq_sud1}{\rm deg}^{\forall\mbox{\rm -}{\rm
sem}}_{\rm rig}(\mathcal{M})\leq{\rm deg}^{\forall\mbox{\rm -}{\rm
synt}}_{\rm rig}(\mathcal{M}),\end{equation} the equality in
(\ref{eq_sud1}) means that either there are no finite sets $A$
with identical $A$-automorphisms only, or minimal finite sets $A$
with identical $A$-automorphisms only have unbounded
cardinalities, or all finite $A\subseteq M$ of some fixed
cardinality $n$ satisfy $M={\rm dcl}(A)$ and some $A$ with $|A|=n$
does not have proper subsets $A'$ such that there are identical
$A'$-automorphisms only;
\begin{equation}\label{eq_sud2}{\rm
deg}^{\exists\mbox{\rm -}{\rm sem}}_{\rm rig}(\mathcal{M})\leq{\rm
deg}^{\exists\mbox{\rm -}{\rm synt}}_{\rm
rig}(\mathcal{M}),\end{equation} the equality in (\ref{eq_sud2})
means that either there are no finite sets $A$ with identical
$A$-automorphisms only, or there is finite $A\subseteq M$ such
that $M={\rm dcl}(A)$, and there are no sets $A'$ with less
cardinalities such that there are identical $A'$-automorphisms
only;
\begin{equation}\label{eq_sud3}{\rm
deg}^{\exists\mbox{\rm -}{\rm sem}}_{\rm rig}(\mathcal{M})\leq{\rm
deg}^{\forall\mbox{\rm -}{\rm sem}}_{\rm
rig}(\mathcal{M}),\end{equation} the equality in (\ref{eq_sud3})
means that either there are no finite sets $A$ with identical
$A$-automorphisms only, or there is finite $A\subseteq M$ with
identical $A$-automorphism only and each finite $A'\subseteq M$
with $|A'|\geq|A|$ has a minimal restriction $A''$, under
inclusion, with $|A''|=|A|$ and with identical $A''$-automorphism
only;
\begin{equation}\label{eq_sud4}{\rm
deg}^{\exists\mbox{\rm -}{\rm synt}}_{\rm
rig}(\mathcal{M})\leq{\rm deg}^{\forall\mbox{\rm -}{\rm
synt}}_{\rm rig}(\mathcal{M}).\end{equation} the equality in
(\ref{eq_sud4}) means that either there are no finite sets $A$
with ${\rm dcl}(A)=M$, or there is finite $A\subseteq M$ with
${\rm dcl}(A)=M$ and each finite $A'\subseteq M$ with
$|A'|\geq|A|$ has a minimal restriction $A''$, under inclusion,
with $|A''|=|A|$ and with ${\rm dcl}(A'')=M$.

\begin{example}\label{ex_rigid_1} {\rm
The structure $\mathcal{M}=\langle\omega,\leq\rangle$ is both
semantically and syntactically rigid, therefore ${\rm
deg}^{\forall\mbox{\rm -}{\rm sem}}_{\rm rig}(\mathcal{M})={\rm
deg}^{\exists\mbox{\rm -}{\rm sem}}_{\rm rig}(\mathcal{M})={\rm
deg}^{\forall\mbox{\rm -}{\rm synt}}_{\rm rig}(\mathcal{M})={\rm
deg}^{\exists\mbox{\rm -}{\rm synt}}_{\rm rig}(\mathcal{M})=0$. We
observe the same effect for arbitrary structures in which each
element is marked by a constant.}
\end{example}

\begin{example}\label{ex_rigid_2} {\rm
If $\mathcal{M}$ has the empty language then $${\rm
deg}^{\forall\mbox{\rm -}{\rm sem}}_{\rm rig}(\mathcal{M})={\rm
deg}^{\exists\mbox{\rm -}{\rm sem}}_{\rm rig}(\mathcal{M})={\rm
deg}^{\forall\mbox{\rm -}{\rm synt}}_{\rm rig}(\mathcal{M})={\rm
deg}^{\exists\mbox{\rm -}{\rm synt}}_{\rm
rig}(\mathcal{M})=|M|-1$$ if $\mathcal{M}$ is finite, and and
these values equal $\infty$ if $\mathcal{M}$ is infinite. }
\end{example}

\begin{example}\label{ex_rigid_21} {\rm
If $\mathcal{V}$ is a vector space over a field $F$ then we have
the following criterion for the semantic/syntactic rigidity: ${\rm
deg}^{\forall\mbox{\rm -}{\rm sem}}_{\rm rig}(\mathcal{V})={\rm
deg}^{\exists\mbox{\rm -}{\rm sem}}_{\rm rig}(\mathcal{V})={\rm
deg}^{\forall\mbox{\rm -}{\rm synt}}_{\rm rig}(\mathcal{V})={\rm
deg}^{\exists\mbox{\rm -}{\rm synt}}_{\rm rig}(\mathcal{V})=0$ iff
${\rm dim}(\mathcal{V})\leq 1$ and $|F|=2$ for ${\rm
dim}(\mathcal{V})=1$. If $\mathcal{V}$ is not rigid then ${\rm
deg}^{\exists\mbox{\rm -}{\rm sem}}_{\rm rig}(\mathcal{V})={\rm
deg}^{\exists\mbox{\rm -}{\rm synt}}_{\rm rig}(\mathcal{V})={\rm
dim}(\mathcal{V})$ for finite ${\rm dim}(\mathcal{V})$, and ${\rm
deg}^{\exists\mbox{\rm -}{\rm sem}}_{\rm rig}(\mathcal{V})={\rm
deg}^{\exists\mbox{\rm -}{\rm synt}}_{\rm
rig}(\mathcal{V})=\infty$, otherwise. Besides, ${\rm
deg}^{\forall\mbox{\rm -}{\rm sem}}_{\rm rig}(\mathcal{V})={\rm
deg}^{\forall\mbox{\rm -}{\rm synt}}_{\rm
rig}(\mathcal{V})=\infty$ if ${\rm dim}(\mathcal{V})$ is infinite,
or ${\rm dim}(\mathcal{V})\geq 1$ and $F$ is infinite. Finally for
${\rm dim}(\mathcal{V})=n\in\omega\setminus\{0\}$ and
$|F|=m\in\omega\setminus\{0\}$ with $(n,m)\ne(1,2)$, we have ${\rm
deg}^{\forall\mbox{\rm -}{\rm sem}}_{\rm rig}(\mathcal{V})={\rm
deg}^{\forall\mbox{\rm -}{\rm synt}}_{\rm
rig}(\mathcal{V})=(n-1)m+1$, since we obtain the rigidity taking
all vectors in a $(n-1)$-dimensional subspace $\mathcal{V}'$, with
$(n-1)m$ elements, and a vector in
$\mathcal{V}\setminus\mathcal{V}'$.}
\end{example}

\begin{example}\label{ex_rigid_3} {\rm
Let $\mathcal{M}$ be a structure of disjoint infinite unary
predicates $P_i$, $i\in I$, expanded by constants for all elements
in $\bigcup\limits_{i\in I}P_i$. Since $\mathcal{M}$ is both
semantically and syntactically rigid we have ${\rm
deg}^{Q\mbox{\rm -}{\rm sem}}_{\rm rig}(\mathcal{M})= {\rm
deg}^{Q\mbox{\rm -}{\rm synt}}_{\rm rig}(\mathcal{M})=0$ for
$Q\in\{\forall,\exists\}$. At the same time extending $n$
predicates $P_i$ by new elements $a_i$ we obtain
$\mathcal{N}\succ\mathcal{M}$ with ${\rm deg}^{\forall\mbox{\rm
-}{\rm sem}}_{\rm rig}(\mathcal{N})={\rm deg}^{\exists\mbox{\rm
-}{\rm sem}}_{\rm rig}(\mathcal{N})=0$, ${\rm
deg}^{\exists\mbox{\rm -}{\rm synt}}_{\rm rig}(\mathcal{N})=n$,
${\rm deg}^{\forall\mbox{\rm -}{\rm synt}}_{\rm
rig}(\mathcal{N})=\infty$. Moreover, if infinitely many $P_i$ are
extended by new elements $a_i$ then the correspondent elementary
extension $\mathcal{N}'$ of $\mathcal{M}$ has the following
characteristics: ${\rm deg}^{\exists\mbox{\rm -}{\rm sem}}_{\rm
rig}(\mathcal{N}')=0$, ${\rm deg}^{\exists\mbox{\rm -}{\rm
synt}}_{\rm rig}(\mathcal{N}')=n$ and ${\rm deg}^{\forall\mbox{\rm
-}{\rm sem}}_{\rm rig}(\mathcal{N}')={\rm deg}^{\forall\mbox{\rm
-}{\rm synt}}_{\rm rig}(\mathcal{N}')=\infty$. Besides, if some
extended $P_i$ are again extended by $m$ new elements in total
then an appropriate elementary extension $\mathcal{N}_{m,n}$ has
the following characteristics: ${\rm deg}^{\exists\mbox{\rm -}{\rm
sem}}_{\rm rig}(\mathcal{N}_{m,n})=m$, ${\rm
deg}^{\exists\mbox{\rm -}{\rm synt}}_{\rm
rig}(\mathcal{N}_{m,n})=m+n$, ${\rm deg}^{\forall\mbox{\rm -}{\rm
sem}}_{\rm rig}(\mathcal{N}_{m,n})={\rm deg}^{\forall\mbox{\rm
-}{\rm synt}}_{\rm rig}(\mathcal{N}_{m,n})=\infty$ including the
possibility ${\rm deg}^{\exists\mbox{\rm -}{\rm sem}}_{\rm
rig}(\mathcal{N}_{\mu,n})={\rm deg}^{\exists\mbox{\rm -}{\rm
synt}}_{\rm rig}(\mathcal{N}_{\mu,n})={\rm deg}^{\forall\mbox{\rm
-}{\rm sem}}_{\rm rig}(\mathcal{N}_{\mu,n})={\rm
deg}^{\forall\mbox{\rm -}{\rm synt}}_{\rm
rig}(\mathcal{N}_{\mu,n})=\infty$ if $\mu\geq\omega$ new elements
are added.}
\end{example}

Thus by Example \ref{ex_rigid_3} the difference between ${\rm
deg}^{\exists\mbox{\rm -}{\rm sem}}_{\rm rig}(\mathcal{M})$ and
${\rm deg}^{\exists\mbox{\rm -}{\rm synt}}_{\rm rig}(\mathcal{M})$
can be arbitrary large. In view of Proposition
\ref{prop_rig_sem_synt} and inequality \ref{eq_sud2} we obtain the
following theorem on distributions for these characteristics:

\begin{theorem}\label{exist_sem_synt}
$1.$ The pairs  $\left({\rm deg}^{\exists\mbox{\rm -}{\rm
sem}}_{\rm rig}(\mathcal{M}), {\rm deg}^{\exists\mbox{\rm -}{\rm
synt}}_{\rm rig}(\mathcal{M})\right)$ belong to the set ${\rm
DEG}^{\rm \exists\mbox{\rm -}{\rm sem}, \exists\mbox{\rm -}{\rm
synt}}_{\rm rig}=\{(\mu,\nu)\mid
\mu,\nu\in\omega\cup\{\infty\},\mu\leq\nu\}$.

$2.$ For each pair $(\mu,\nu)\in{\rm DEG}^{\rm \exists\mbox{\rm
-}{\rm sem}, \exists\mbox{\rm -}{\rm synt}}_{\rm rig}$ there
exists a structure $\mathcal{M}_{\mu,\nu}$ such that $${\rm
deg}^{\exists\mbox{\rm -}{\rm sem}}_{\rm
rig}(\mathcal{M}_{\mu,\nu})=\mu, \, {\rm deg}^{\exists\mbox{\rm
-}{\rm synt}}_{\rm rig}(\mathcal{M}_{\mu,\nu})=\nu.$$
\end{theorem}

Example \ref{ex_rigid_3} shows that values in ${\rm DEG}^{\rm
\exists\mbox{\rm -}{\rm sem}, \exists\mbox{\rm -}{\rm synt}}_{\rm
rig}$ in Theorem \ref{exist_sem_synt} are covered by structures in
countable languages $\Sigma_1$ of unary predicates. Now we
describe possibilities for the pairs $\left({\rm
deg}^{\forall\mbox{\rm -}{\rm sem}}_{\rm rig}(\mathcal{M}), {\rm
deg}^{\forall\mbox{\rm -}{\rm synt}}_{\rm
rig}(\mathcal{M})\right)$ in these languages $\Sigma_1$.

\begin{proposition}\label{forall_sem_synt1}
For any structure $\mathcal{M}$ in a language $\Sigma_1$ of unary
predicates the pair
$$\left({\rm deg}^{\forall\mbox{\rm -}{\rm sem}}_{\rm
rig}(\mathcal{M}), {\rm deg}^{\forall\mbox{\rm -}{\rm synt}}_{\rm
rig}(\mathcal{M})\right)$$ has one of the following possibilities:

$1)$ $(0,0)$, if $\mathcal{M}$ is both semantically and
syntactically rigid;

$2)$ $(n,n)$, if $\mathcal{M}$ is finite with $n+1$ elements and
it is not semantically rigid that is not syntactically rigid;

$3)$ $(0,\infty)$, if $\mathcal{M}$ is infinite, semantically
rigid but not syntactically rigid;

$4)$ $(\infty,\infty)$, if $\mathcal{M}$ is infinite and both not
semantically rigid and not syntactically rigid.
\end{proposition}

Proof. If $\mathcal{M}$ is syntactically rigid then we have
$\left({\rm deg}^{\forall\mbox{\rm -}{\rm sem}}_{\rm
rig}(\mathcal{M}), {\rm deg}^{\forall\mbox{\rm -}{\rm synt}}_{\rm
rig}(\mathcal{M})\right)=(0,0)$ by the inequality (\ref{eq_sud1}).
Now we assume that $\mathcal{M}$ is not syntactically rigid and
consider the following cases.

Case 1: $\mathcal{M}$ is semantically rigid, i.e., ${\rm
deg}^{\forall\mbox{\rm -}{\rm sem}}_{\rm rig}(\mathcal{M})=0$. In
such a case $\mathcal{M}$ is infinite since finite structures have
isolated $1$-types only and there are complete $1$-types over
empty set with at least two realizations that contradicts the
semantic rigidity for the language $\Sigma_1$. Again using the
unary language $\Sigma_1$ and the arguments of \cite[Section
8.1]{ErPa} that all $1$-types, over empty set, are forced by
formulae of quantifier free diagrams and formulae describing
estimations for cardinalities of their solutions, with independent
actions of automorphisms in distinct sets of realizations of
$1$-types. Thus each $1$-type has at most one realization in
$\mathcal{M}$. Since $\mathcal{M}$ is not syntactically rigid,
$\mathcal{M}$ realizes at least one nonisolated $1$-type $p(x)$ by
some unique element $a$. Now for any $n\in\omega$ we can take $n$
realizations of other $1$-types forming a set $A$ such that
$a\notin{\rm dcl}(A)$. It implies ${\rm deg}^{\forall\mbox{\rm
-}{\rm synt}}_{\rm rig}(\mathcal{M})=\infty$.

Case 2: $\mathcal{M}$ is not semantically rigid and
$|M|=n+1\in\omega$. In such a case $\mathcal{M}$ has a complete
$1$-type $p(x)$ with at least two realizations $a$ and $b$. Since
there is an $(M\setminus\{a,b\})$-automorphism $f$ with $f(a)=b$,
we obtain ${\rm deg}^{\forall\mbox{\rm -}{\rm sem}}_{\rm
rig}(\mathcal{M})=n$ implying ${\rm deg}^{\forall\mbox{\rm -}{\rm
synt}}_{\rm rig}(\mathcal{M})=n$ by the inequality (\ref{eq_sud1})
and the syntactic rigidity of $\mathcal{M}$ over each $n$-element
set.

Case 3: $\mathcal{M}$ is not semantically rigid and it is
infinite. In such a case $\mathcal{M}$ has a complete $1$-type
$p(x)$ with at least two realizations $a$ and $b$ and such that
realizations of other $1$-types allow to form arbitrarily large
finite set $A$ such that some $A$-automorphism transforms $a$ in
$b$. It means that ${\rm deg}^{\forall\mbox{\rm -}{\rm sem}}_{\rm
rig}(\mathcal{M})=\infty$ implying ${\rm deg}^{\forall\mbox{\rm
-}{\rm synt}}_{\rm rig}(\mathcal{M})=\infty$ by the inequality
(\ref{eq_sud1}).
\endproof

\medskip
Combining arguments for Theorems \ref{exist_sem_synt} and
\ref{forall_sem_synt1} we obtain the following possibilities for
tetrads ${\rm deg}_4(\mathcal{M})\rightleftharpoons\left({\rm
deg}^{\exists\mbox{\rm -}{\rm sem}}_{\rm rig}(\mathcal{M}), {\rm
deg}^{\exists\mbox{\rm -}{\rm synt}}_{\rm rig}(\mathcal{M}),{\rm
deg}^{\forall\mbox{\rm -}{\rm sem}}_{\rm rig}(\mathcal{M}), {\rm
deg}^{\forall\mbox{\rm -}{\rm synt}}_{\rm
rig}(\mathcal{M})\right)$ in a language of unary predicates:

\begin{corollary}\label{cor_sem_synt1}
For any structure $\mathcal{M}$ in a language $\Sigma_1$ of unary
predicates the tetrad ${\rm deg}_4(\mathcal{M})$ has one of the
following possibilities:

$1)$ $(0,0,0,0)$, if $\mathcal{M}$ is both semantically and
syntactically rigid;

$2)$ $(m,m,n,n)$, if $\mathcal{M}$ is finite with $n+1$ elements
and it is not semantically rigid that is not syntactically rigid
with some minimal $m$-elements set $A\subset M$, $1\leq m\leq n$,
producing ${\rm dcl}(A)=M$;

$3)$ $(0,\nu,0,\infty)$, if $\mathcal{M}$ is infinite,
semantically rigid but not syntactically rigid, with
$1\leq\nu\leq\infty$;

$4)$ $(\mu,\nu,\infty,\infty)$, if $\mathcal{M}$ is infinite and
both not semantically rigid and not syntactically rigid, with
$1\leq\mu\leq\nu\leq\infty$.
\end{corollary}

\begin{example}\label{ex_rigid_4} {\rm
Let $\mathcal{M}$ be a finitely generated algebra by a set $X$.
Then by the definition we have ${\rm deg}^{\exists\mbox{\rm -}{\rm
synt}}_{\rm rig}(\mathcal{M})\leq |X|$ which implies ${\rm
deg}^{\exists\mbox{\rm -}{\rm sem}}_{\rm rig}(\mathcal{M})\leq
|X|$ by the inequality (\ref{eq_sud2}). Here, if additionally the
generating set $X$ admits substitutions by any $Y\subseteq M$ with
$|Y|=|X|$ and these substitutions preserve the generating property
then we have ${\rm deg}^{\forall\mbox{\rm -}{\rm synt}}_{\rm
rig}(\mathcal{M})\leq |X|$ which implies ${\rm
deg}^{\exists\mbox{\rm -}{\rm sem}}_{\rm rig}(\mathcal{M})\leq
|X|$ by the inequality (\ref{eq_sud1}). For instance, if
$\mathcal{M}$ is a directed graph forming a finite cycle of
positive length then ${\rm deg}_4(\mathcal{M})=(1,1,1,1)$.

Since algebras, with constants and unary operations, can define
arbitrary configurations of unary predicates, possibilities for
characteristics ${\rm deg}_4(\mathcal{M})$ in Corollary
\ref{cor_sem_synt1} can be realized in the class of algebras,
too.}
\end{example}

\begin{example}\label{ex_rigid_5} {\rm Let ${\rm pm}={\rm pm}(G_1,G_2,\mathcal{P})$ be a connected
polygonometry of a group pair $(G_1,G_2)$ on an exact pseudoplane
$\mathcal{P}$, and $\mathcal{M}=\mathcal{M}({\rm pm})$ be a
ternary structure for ${\rm pm}$ \cite{GP}. Since all points $a$
in $\mathcal{M}$ are connected by automorphisms we have ${\rm
acl}(\{a\})=\{a\}$. At the same time any two distinct points
$a,b\in M({\rm pm})$ (laying in a common line) define all points
in $\mathcal{M}$ by line and angle parameters of broken lines. It
implies $M({\rm pm})={\rm dcl}(\{a,b\})$. If line and angle
parameters of shortest broken lines connecting arbitrary distinct
points $a$ and $b$ are defined uniquely then $M({\rm pm})={\rm
dcl}(\{a,b\})$ for these points, too.  Hence, in such a case,
${\rm deg}^{\exists\mbox{\rm -}{\rm sem}}_{\rm
rig}(\mathcal{M})={\rm deg}^{\exists\mbox{\rm -}{\rm synt}}_{\rm
rig}(\mathcal{M})={\rm deg}^{\forall\mbox{\rm -}{\rm sem}}_{\rm
rig}(\mathcal{M})={\rm deg}^{\forall\mbox{\rm -}{\rm synt}}_{\rm
rig}(\mathcal{M})\leq 2$. Moreover, these degree values equal $1$
iff ${\rm pm}$ consists of unique line and with at least two
points, i.e., $|G_1|>1$ and $|G_2|=1$. Finally, for a
polygonometry ${\rm pm}$, the degrees equal $0$ iff ${\rm pm}$
consists of unique point.

If parameters of broken lines do not define these broken lines by
endpoints then finite cardinalities of points in these lines can
be unbounded. Indeed, taking opposite vertices $a$ and $b$ in an
$n$-cube \cite{GP, mct} or in its polygonometry ${\rm pm}$ we
obtain $n$ adjacent vertices $c_1,\ldots,c_n$ for $a$ and these
vertices are connected by $\{a,b\}$-automorphisms. Moreover, in
such a case, ${\rm deg}^{\exists\mbox{\rm -}{\rm sem}}_{\rm
rig}(\mathcal{M})={\rm deg}^{\exists\mbox{\rm -}{\rm synt}}_{\rm
rig}(\mathcal{M})=n+1$ witnessed, for instance, by the set
$A=\{a,b,c_1,\ldots,c_{n-1}\}$.

The value ${\rm deg}_4(\mathcal{M}_2)=(2,2,2,2)$ for
$\mathcal{M}_2=\mathcal{M}({\rm pm})$ can be increased till ${\rm
deg}_4(\mathcal{M}_n)=(n,n,n,n)$, $n\geq 3$, generalizing group
trigonometries in the following way. We construct a
$(n+1)$-dimensional space consisting of points and $n$-dimensional
hyperplanes. We introduce an incidence $n$-ary relation $I_n$ for
$n$ distinct points to lay on a common hyperplane. Now fixing a
hyperplane $H$ and $n-1$ pairwise distinct points
$a_1,\ldots,a_{n-1}\in H$ we define an exact transitive action of
a group $G_1$ on $H\setminus\{a_1,\ldots,a_{n-1}\}$, i.e., on $H$
with respect to $a_1,\ldots,a_{n-1}$, such that this action is
transformed for any pairwise distinct points
$a'_1,\ldots,a'_{n-1}\in H$. Since each $H$ can be defined by its
$n-1$ distinct points with actions, we can fix
$a_1,\ldots,a_{n-1}$ and move $a_n\in
H\setminus\{a_1,\ldots,a_{n-1}\}$ into points $a'_n$ in other
hyperplanes $H'$ containing $a_1,\ldots,a_{n-1}$. Collecting these
movements we define an action of a group $G_2$ on that bundle of
hyperplanes containing $a_1,\ldots,a_{n-1}$. Then we spread
actions of $G_1$ and $G_2$ for any hyperplanes and bundles of
hyperplanes, respectively, such that all pairwise distinct
$a_1,\ldots,a_{n-1}$ and $a'_1,\ldots,a'_{n-1}$ are connected by
automorphisms with respect to these actions.

For instance, taking the set $P$ of planes in $\mathbb R^3$, a
plane $\pi\in P$ and distinct points $a_1,a_2\in P$ the action of
$G_1$ can be defined as $\mathbb R\times A$ with the side group
$\mathbb R$ and angle group $A$ defining both the directed
distance $d\in\mathbb R$ from $a_1$ to a point $a_3\in\pi$ and the
angle value $\alpha$ from the side $a_1\hat{\,}a_2$ to the side
$a_1\hat{\,}a_3$. And $G_2$ is the rotation group for the planes
in $P$ around the lines $l(a_1,a_2)$.

Now we extend the language $\{I_n\}$ by $(n+1)$-ary predicates
$Q_{g_1}$, $g_1\in G_1$, such that first $(n-1)$-coordinates
$\overline{a}$ in $\langle\overline{a},b,c\rangle\in Q_{g_1}$ are
exhausted by $a_1,\ldots,a_{n-1}$ and $c=bg_1$ with respect to
$a_1,\ldots,a_{n-1}$. Simultaneously we define predicates
$R_{g_2}$, $g_2\in G_2$, of arities $n+1$ such that each $R_{g_2}$
realizes a rotation of a hyperplane with respect to
$a_1,\ldots,a_{n-1}$ by the element $g_2$. We obtain a structure
$\mathcal{M}_n$ whose values ${\rm deg}^{Q\mbox{\rm -}{\rm
sem}}_{\rm rig}(\mathcal{M}_n)$ and ${\rm deg}^{Q\mbox{\rm -}{\rm
synt}}_{\rm rig}(\mathcal{M}_n)$, for $Q\in\{\forall,\exists\}$
equal $n$.

The construction above admits a generalization for polygonometries
${\rm pm}(G_1,G_2,\mathcal{P})$ of group pairs transforming
$(G_1,G_2)$ a pseudoplane $\mathcal{P}$ to a pseudospace
$\mathcal{S}$ with hyperplanes $H$ such that $H={\rm
dcl}(\{a_1,\ldots,a_{n}\})$ for any pairwise distinct points
$a_1,\ldots,a_{n}\in H$ and with ${\rm
dcl}(\{b_1,\ldots,b_{n-1}\})=\{b_1,\ldots,b_{n-1}\}$ for any
$b_1,\ldots,b_{n-1}\in\mathcal{S}$.}
\end{example}

Comparing characteristics ${\rm deg}^{\exists\mbox{\rm -}{\rm
sem}}_{\rm rig}(\mathcal{M})$ / ${\rm deg}^{\exists\mbox{\rm
-}{\rm synt}}_{\rm rig}(\mathcal{M})$ and ${\rm
deg}^{\forall\mbox{\rm -}{\rm sem}}_{\rm rig}(\mathcal{M})$ /
${\rm deg}^{\forall\mbox{\rm -}{\rm synt}}_{\rm rig}(\mathcal{M})$
we observe that the first ones produce cardinalities of ``best'',
i.e., minimal sets generating the structure $\mathcal{M}$ and the
second ones give cardinalities of ``worst'' generating sets. It is
natural to describe possibilities of ``intermediate'' generating
sets. For this aim we define the degrees of rigidity with respect
to a subset $A$ of $M$ as follows:

\medskip
{\bf Definition.} For a set $A$ in $\mathcal{M}$ and an expansion
$\mathcal{M}_A$ of  $\mathcal{M}$ by constants in $A$, the least
$n$ such that $\mathcal{M}_A$ is $Q$-semantically /
$Q$-syntactically $n$-rigid, where $Q\in\{\forall,\exists\}$, is
called the {\em $(Q,A)$-semantical / $(Q,A)$-syntactical degree of
rigidity}, it is denoted by ${\rm deg}^{Q\mbox{\rm -}{\rm
sem}}_{{\rm rig},A}(\mathcal{M})$ and ${\rm deg}^{Q\mbox{\rm
-}{\rm synt}}_{{\rm rig},A}(\mathcal{M})$, respectively. If such
$n$ does not exists we put ${\rm deg}^{Q\mbox{\rm -}{\rm
sem}}_{{\rm rig},A}(\mathcal{M})=\infty$ and ${\rm
deg}^{Q\mbox{\rm -}{\rm synt}}_{{\rm rig},A}(\mathcal{M})=\infty$,
respectively.

Any expansion $\mathcal{M}_A$ of $\mathcal{M}$ with ${\rm
deg}^{\exists\mbox{\rm -}s}_{\rm rig}(\mathcal{M}_A)=0$, for
$s\in\{{\rm sem},{\rm synt}\}$, is called a {\em
$s$-rigiditization} or simply a {\em rigiditization} of
$\mathcal{M}$.

\medskip
We have the following properties for $(Q,A)$-semantical and
$(Q,A)$-syntactical degrees of rigidity:

\begin{proposition}\label{deg_A} Let $\mathcal{M}$ be a structure, $A\subseteq M$, $Q\in\{\forall,\exists\}$, $s\in\{{\rm sem},{\rm synt}\}$.
Then the following assertions hold:

$1.$ {\rm (Preservation of degrees of rigidity)} If
$A\subseteq{\rm dcl}(\emptyset)$ then ${\rm deg}^{Q\mbox{\rm
-}s}_{\rm rig}(\mathcal{M})={\rm deg}^{Q\mbox{\rm -}s}_{{\rm
rig},A}(\mathcal{M})$.

$2.$ {\rm (Rigiditization)} If $A$ contains a witnessing set for
the finite value ${\rm deg}^{\exists\mbox{\rm -}s}_{\rm
rig}(\mathcal{M})$ then ${\rm deg}^{\exists\mbox{\rm -}s}_{{\rm
rig},A}(\mathcal{M})=0$.

$3.$ {\rm (Monotony)} If $A\subseteq B\subseteq M$ then ${\rm
deg}^{Q\mbox{\rm -}s}_{{\rm rig},A}(\mathcal{M})\geq{\rm
deg}^{Q\mbox{\rm -}s}_{{\rm rig},B}(\mathcal{M})$.

$4.$ {\rm (Additivity)} If $A$ witnesses the finite value ${\rm
deg}^{\exists\mbox{\rm -}s}_{\rm rig}(\mathcal{M})$ then for any
$A'\subseteq A$, $${\rm deg}^{\exists\mbox{\rm -}s}_{\rm
rig}(\mathcal{M})={\rm deg}^{\exists\mbox{\rm -}s}_{{\rm
rig},A'}(\mathcal{M})+{\rm deg}^{\exists\mbox{\rm -}s}_{{\rm
rig},A\setminus A'}(\mathcal{M}).$$

$5.$ {\rm (Cofinite character)} If $A$ is cofinite in
$\mathcal{M}$ then ${\rm deg}^{\exists\mbox{\rm -}{\rm sem}}_{{\rm
rig},A}(\mathcal{M})$ and ${\rm deg}^{\exists\mbox{\rm -}{\rm
synt}}_{{\rm rig},A}(\mathcal{M})$ are natural.

$6.$ {\rm (Finite rigiditization)} Any cofinite set $A$ in
$\mathcal{M}$ has a minimal finite extension $A'$ such that
$\mathcal{M}_{A'}$ is semantically / syntactically rigid.

\end{proposition}

Proof. 1. If $A\subseteq{\rm dcl}(\emptyset)$ then ${\rm
Aut}(\mathcal{M})={\rm Aut}(\mathcal{M}_A)$ and therefore the
equalities ${\rm deg}^{Q\mbox{\rm -}s}_{\rm rig}(\mathcal{M})={\rm
deg}^{Q\mbox{\rm -}s}_{{\rm rig},A}(\mathcal{M})$ hold for $s={\rm
sem}$. For the case $s={\rm synt}$ the required equalities are
satisfied in view of ${\rm dcl}(B)={\rm dcl}(A\cup B)$ for any
$B\subseteq M$.

2. If $A$ contains a witnessing set for the finite value ${\rm
deg}^{\exists\mbox{\rm -}{\rm sem}}_{\rm rig}(\mathcal{M})$ then
there exists identical $A$-automorphism of $\mathcal{M}$ only
implying ${\rm deg}^{\exists\mbox{\rm -}{\rm sem}}_{{\rm
rig},A}(\mathcal{M})=0$. Similarly if $A$ contains a witnessing
set for the finite value ${\rm deg}^{\exists\mbox{\rm -}{\rm
synt}}_{\rm rig}(\mathcal{M})$ then ${\rm dcl}(A)=M$ producing
${\rm deg}^{\exists\mbox{\rm -}{\rm synt}}_{{\rm
rig},A}(\mathcal{M})=0$.

3. If $A\subseteq B\subseteq M$ then ${\rm
Aut}(\mathcal{M}_B)\leq{\rm Aut}(\mathcal{M}_A)$ therefore the
inequalities ${\rm deg}^{Q\mbox{\rm -}s}_{{\rm
rig},A}(\mathcal{M})\geq{\rm deg}^{Q\mbox{\rm -}s}_{{\rm
rig},B}(\mathcal{M})$ hold for $s={\rm sem}$. For the case $s={\rm
synt}$ the required equalities are satisfied in view of ${\rm
dcl}(A\cup C)\subseteq{\rm dcl}(B\cup C)$ for any $C\subseteq M$.

4. If $A$ witnesses the finite value ${\rm deg}^{\exists\mbox{\rm
-}s}_{\rm rig}(\mathcal{M})$ then we divide $A$ into two disjoint
parts $A_1$ and $A_2$ and by the definition of ${\rm
deg}^{\exists\mbox{\rm -}s}_{\rm rig}(\mathcal{M})$, both $A_1$
and $A_2$ are extended till minimal $A$ witnessing the semantic /
syntactic rigidity. Thus $A_1$ witnesses the value ${\rm
deg}^{\exists\mbox{\rm -}{\rm sem}}_{\rm rig}(\mathcal{M}_{A_2})$
and $A_2$ witnesses the value ${\rm deg}^{\exists\mbox{\rm -}{\rm
sem}}_{\rm rig}(\mathcal{M}_{A_1})$ producing the required
equation ${\rm deg}^{\exists\mbox{\rm -}s}_{\rm
rig}(\mathcal{M})={\rm deg}^{\exists\mbox{\rm -}s}_{{\rm
rig},A'}(\mathcal{M})+{\rm deg}^{\exists\mbox{\rm -}s}_{{\rm
rig},A\setminus A'}(\mathcal{M}).$

5. If $A$ is cofinite in $\mathcal{M}$ then there are only
finitely many elements, all in $M\setminus A$, witnessing the
values ${\rm deg}^{\exists\mbox{\rm -}{\rm sem}}_{{\rm
rig},A}(\mathcal{M})$ and ${\rm deg}^{\exists\mbox{\rm -}{\rm
synt}}_{{\rm rig},A}(\mathcal{M})$. Thus these values are natural.

6. It is immediately implied by Items 2 and 5.
\endproof

\medskip
In view of Proposition \label{deg_A} fixing a subset in
$\mathcal{M}$ large enough we obtain its rigiditization. At the
same time the following assertion clarifies that small subsets can
produce the rigiditization for structures in bounded cardinalities
only.

\begin{proposition}\label{card_mod}
$1.$ If ${\rm deg}^{\exists\mbox{\rm -}{\rm synt}}_{\rm
rig}(\mathcal{M})$ is finite then $|M|\leq{\rm
max}\{\Sigma(\mathcal{M}),\omega\}$.

$1.$ If $\mathcal{M}$ is homogeneous and ${\rm
deg}^{\exists\mbox{\rm -}{\rm sem}}_{\rm rig}(\mathcal{M})$ is
finite then $|M|\leq 2^{{\rm max}\{\Sigma(\mathcal{M}),\omega\}}$.
\end{proposition}

Proof. 1. If ${\rm deg}^{\exists\mbox{\rm -}{\rm synt}}_{\rm
rig}(\mathcal{M})$ is finite then there is a finite set
$A\subseteq M$ witnessing that value, with $M={\rm dcl}(A)$. This
equality is witnessed by at most by ${\rm
max}\{\Sigma(\mathcal{M}),\omega\}$ formulae such that each
element in $\mathcal{M}$ is defined by a formula in the language
$\Sigma(\mathcal{M}_A)$. Since there are ${\rm
max}\{\Sigma(\mathcal{M}),\omega\}$
$\Sigma(\mathcal{M}_A)$-formulae we obtain at most ${\rm
max}\{\Sigma(\mathcal{M}),\omega\}$ elements in $\mathcal{M}$.

2. If a finite set $A\subseteq M$ witnesses the finite value ${\rm
deg}^{\exists\mbox{\rm -}{\rm sem}}_{\rm rig}(\mathcal{M})$ and
$\mathcal{M}$ is homogeneous possibilities for $A$-automorphisms
fixing elements of $\mathcal{M}$ are exhausted by single
realizations of types in $S^1(A)$. Since there are at most
$2^{{\rm max}\{\Sigma(\mathcal{M}),\omega\}}$ these types that
value is the required upper bound for the cardinality of
semantically rigid structure $\mathcal{M}_A$.
\endproof

\medskip
Proposition \ref{card_mod} immediately implies the following:

\begin{corollary}\label{cor_card_mod}
$1.$ If ${\rm deg}^{\exists\mbox{\rm -}{\rm synt}}_{{\rm
rig},A}(\mathcal{M})$ is finite then $|M|\leq{\rm
max}\{\Sigma(\mathcal{M}),|A|,\omega\}$.

$1.$ If $\mathcal{M}$ is homogeneous and ${\rm
deg}^{\exists\mbox{\rm -}{\rm sem}}_{{\rm rig},A}(\mathcal{M})$ is
finite then $|M|\leq 2^{{\rm
max}\{\Sigma(\mathcal{M}),|A|,\omega\}}$.
\end{corollary}

\section{Indexes of rigidity}

{\bf Definition.} For a set $A$ in a structure $\mathcal{M}$ the
{\em index of rigidity} of $\mathcal{M}$ over $A$, denoted by
${\rm ind}_{\rm rig}(\mathcal{M}/A)$ is the supremum of
cardinalities for the set of solutions of algebraic types ${\rm
tp}(a/A)$ for $a\in M$. We put ${\rm ind}_{\rm
rig}(\mathcal{M})={\rm ind}_{\rm rig}(\mathcal{M}/\emptyset)$.
Here we assume that ${\rm ind}_{\rm rig}(\mathcal{M})=0$ if
$\mathcal{M}$ does not have algebraic types ${\rm tp}(a)$ for
$a\in M$.

\begin{note}\label{rem_ind}{\rm
By the definition we have ${\rm ind}_{\rm
rig}(\mathcal{M}/A)\in\omega+1$.}
\end{note}

\begin{example}\label{ind_rig} {\rm
1. If $\mathcal{M}$ is a structure of unary predicates $P_i$,
$i\in I$, then ${\rm ind}_{\rm rig}(\mathcal{M})=0$ iff there are
no finite nonempty intersections $P^{\delta_1}_{i_1}\cap\ldots\cap
P^{\delta_k}_{i_k}$, $\delta_1,\ldots,\delta_k\in\{0,1\}$. We have
${\rm ind}_{\rm rig}(\mathcal{M})=1$ iff ${\rm
dcl}(\emptyset)\ne\emptyset$ and there are no maximal finite
intersections $P^{\delta_1}_{i_1}\cap\ldots\cap
P^{\delta_k}_{i_k}$ with at least two elements. Besides, ${\rm
ind}_{\rm rig}(\mathcal{M})\in\omega$ iff these finite
intersections have bounded cardinalities, and all natural
possibilities $n$ are realized by predicates with exactly $n$
elements and infinite complements. Otherwise, i.e., for ${\rm
ind}_{\rm rig}(\mathcal{M})=\omega$, these finite intersections
have unbounded cardinalities.

2. If $\mathcal{M}$ is a structure of an equivalence relation $E$,
then ${\rm ind}_{\rm rig}(\mathcal{M})=0$ iff there are no finite
$E$-classes. We have ${\rm ind}_{\rm rig}(\mathcal{M})=1$ iff
${\rm dcl}(\emptyset)\ne\emptyset$ and there are no finite
$E$-classes with at least two elements. Besides, ${\rm ind}_{\rm
rig}(\mathcal{M})\in\omega$ iff these $E$-classes have bounded
cardinalities, and all natural possibilities $n$ are realized by
infinitely many $E$-classes with exactly $n$ elements. Otherwise,
i.e., for ${\rm ind}_{\rm rig}(\mathcal{M})=\omega$, these
$E$-classes have unbounded cardinalities.

3. If $\mathcal{M}=\mathcal{M}({\rm pm})$ for a polygonometry
${\rm pm}$ then ${\rm ind}_{\rm rig}(\mathcal{M})=0$ iff ${\rm
pm}$ has infinitely many points. Otherwise, if ${\rm pm}$ has
$n\in\omega$ points then ${\rm ind}_{\rm rig}(\mathcal{M})=n$.

More generally, we have the following possibilities for a model
$\mathcal{M}$ of transitive theory $T$, i.e., of a theory with
$|S^1(\emptyset)|=1$:

i) ${\rm ind}_{\rm rig}(\mathcal{M})=0$, if $\mathcal{M}$ is
infinite;

ii) ${\rm ind}_{\rm rig}(\mathcal{M})=|\mathcal{M}|$, if
$\mathcal{M}$ is finite. }
\end{example}

In view of Remark \ref{rem_ind} the following assertion describes
possibilities of indexes of rigidity:

\begin{proposition}\label{th_ind_rig}
For any $\lambda\in\omega+1$ there is a structure
$\mathcal{M}_\lambda$ such that ${\rm ind}_{\rm
rig}(\mathcal{M}_\lambda)=\lambda$.
\end{proposition}

Proof follows by Example \ref{ind_rig}. \endproof

\section{Variations of rigidity for disjoint unions of structures}

{\bf Definition} \cite{Wo}. The {\em disjoint
union}\index{Disjoint union!of structures}
$\bigsqcup\limits_{n\in\omega}\mathcal{
M}_n$\index{$\bigsqcup\limits_{n\in\omega}\mathcal{ M}_n$} of
pairwise disjoint structures $\mathcal{ M}_n$ for pairwise
disjoint predicate languages $\Sigma_n$, $n\in\omega$, is the
structure of language
$\bigcup\limits_{n\in\omega}\Sigma_n\cup\{P^{(1)}_n\mid
n\in\omega\}$ with the universe $\bigsqcup\limits_{n\in\omega}
M_n$, $P_n=M_n$, and interpretations of predicate symbols in
$\Sigma_n$ coinciding with their interpretations in $\mathcal{
M}_n$, $n\in\omega$. The {\em disjoint union of
theories}\index{Disjoint union!of theories} $T_n$ for pairwise
disjoint languages $\Sigma_n$ accordingly, $n\in\omega$, is the
theory
$$\bigsqcup\limits_{n\in\omega}T_n\rightleftharpoons{\rm Th}\left(\bigsqcup\limits_{n\in\omega}\mathcal{ M}_n\right),$$
where\index{$\bigsqcup\limits_{n\in\omega}T_n$} $\mathcal{
M}_n\models T_n$, $n\in\omega$.

\begin{theorem}\label{th_dis}
For any disjoint predicate structures $\mathcal{ M}_1$ and
$\mathcal{ M}_2$, and $s\in\{{\rm sem},{\rm synt}\}$ the following
conditions hold:

$1.$ ${\rm deg}^{\exists\mbox{\rm -}s}_{\rm
rig}(\mathcal{M}_1\sqcup\mathcal{M}_2)={\rm deg}^{\exists\mbox{\rm
-}s}_{\rm rig}(\mathcal{M}_1)+{\rm deg}^{\exists\mbox{\rm
-}s}_{\rm rig}(\mathcal{M}_2)$, in particular, ${\rm
deg}^{\exists\mbox{\rm -}s}_{\rm
rig}(\mathcal{M}_1\sqcup\mathcal{M}_2)$ is finite iff ${\rm
deg}^{\exists\mbox{\rm -}s}_{\rm rig}(\mathcal{M}_1)$ and ${\rm
deg}^{\exists\mbox{\rm -}s}_{\rm rig}(\mathcal{M}_2)$ are finite.

$2.$ ${\rm deg}^{\forall\mbox{\rm -}s}_{\rm
rig}(\mathcal{M}_1\sqcup\mathcal{M}_2)=0$ iff ${\rm
deg}^{\forall\mbox{\rm -}s}_{\rm rig}(\mathcal{M}_1)=0$ and ${\rm
deg}^{\forall\mbox{\rm -}s}_{\rm rig}(\mathcal{M}_2)=0$.

$3.$ If ${\rm deg}^{\forall\mbox{\rm -}s}_{\rm
rig}(\mathcal{M}_1\sqcup\mathcal{M}_2)>0$ then it is finite iff
${\rm deg}^{\forall\mbox{\rm -}s}_{\rm rig}(\mathcal{M}_1)>0$ is
finite and $\mathcal{M}_2$ is finite, or ${\rm
deg}^{\forall\mbox{\rm -}s}_{\rm rig}(\mathcal{M}_2)>0$ is finite
and $\mathcal{M}_1$ is finite. Here,
$${\rm deg}^{\forall\mbox{\rm -}s}_{\rm
rig}(\mathcal{M}_1\sqcup\mathcal{M}_2)={\rm max}\{|M_1|+{\rm
deg}^{\forall\mbox{\rm -}s}_{\rm rig}(\mathcal{M}_2),|M_2|+{\rm
deg}^{\forall\mbox{\rm -}s}_{\rm rig}(\mathcal{M}_1)\}.$$

\end{theorem}

Proof. 1. Let $A_i\subset M_i$ be sets witnessing values ${\rm
deg}^{\exists\mbox{\rm -}s}_{\rm rig}(\mathcal{M}_i)$, $i=1,2$. By
the definition of $\mathcal{M}_1\sqcup\mathcal{M}_2$, $A_1$ and
$A_2$ are disjoint and $A_1\cup A_2$ witnesses the value ${\rm
deg}^{\exists\mbox{\rm -}s}_{\rm
rig}(\mathcal{M}_1\sqcup\mathcal{M}_2)$. Thus ${\rm
deg}^{\exists\mbox{\rm -}s}_{\rm
rig}(\mathcal{M}_1\sqcup\mathcal{M}_2)={\rm deg}^{\exists\mbox{\rm
-}s}_{\rm rig}(\mathcal{M}_1)+{\rm deg}^{\exists\mbox{\rm
-}s}_{\rm rig}(\mathcal{M}_2)$.

2. If ${\rm deg}^{\forall\mbox{\rm -}s}_{\rm
rig}(\mathcal{M}_1\sqcup\mathcal{M}_2)=0$ then the empty set
witnesses that $\mathcal{M}_1\sqcup\mathcal{M}_2$, $\mathcal{M}_1$
and $\mathcal{M}_2$ are $s$-rigid, i.e., rigid with respect to
$s$, implying ${\rm deg}^{\forall\mbox{\rm -}s}_{\rm
rig}(\mathcal{M}_1)=0$ and ${\rm deg}^{\forall\mbox{\rm -}s}_{\rm
rig}(\mathcal{M}_2)=0$. Conversely, if ${\rm
deg}^{\forall\mbox{\rm -}s}_{\rm rig}(\mathcal{M}_1)=0$ and ${\rm
deg}^{\forall\mbox{\rm -}s}_{\rm rig}(\mathcal{M}_2)=0$ then the
empty set witnesses that $\mathcal{M}_1$ and $\mathcal{M}_2$ are
$s$-rigid. Now by the definition of
$\mathcal{M}_1\sqcup\mathcal{M}_2$ we observe that
$\mathcal{M}_1\sqcup\mathcal{M}_2$ is $s$-rigid, too, implying
${\rm deg}^{\forall\mbox{\rm -}s}_{\rm
rig}(\mathcal{M}_1\sqcup\mathcal{M}_2)=0$.

3. Let ${\rm deg}^{\forall\mbox{\rm -}s}_{\rm
rig}(\mathcal{M}_1\sqcup\mathcal{M}_2)>0$ be finite, then by Item
2, ${\rm deg}^{\forall\mbox{\rm -}s}_{\rm rig}(\mathcal{M}_1)>0$
or ${\rm deg}^{\forall\mbox{\rm -}s}_{\rm rig}(\mathcal{M}_2)>0$.
Assuming that ${\rm deg}^{\forall\mbox{\rm -}s}_{\rm
rig}(\mathcal{M}_i)>0$ we can not witness that value by subsets of
$M_{3-i}$, $i=1,2$. Thus $M_{3-i}$ should be finite. Conversely,
let ${\rm deg}^{\forall\mbox{\rm -}s}_{\rm rig}(\mathcal{M}_1)>0$
be finite and $\mathcal{M}_2$ be finite, or ${\rm
deg}^{\forall\mbox{\rm -}s}_{\rm rig}(\mathcal{M}_2)>0$ be finite
and $\mathcal{M}_1$ be finite. Then we can take ${\rm
deg}^{\forall\mbox{\rm -}s}_{\rm rig}(\mathcal{M}_1)$ elements of
$M_1$ and all elements of $M_2$ obtaining the $s$-rigidity of
$\mathcal{M}_1\sqcup\mathcal{M}_2$. Similarly we can take ${\rm
deg}^{\forall\mbox{\rm -}s}_{\rm rig}(\mathcal{M}_2)$ elements of
$M_2$ and all elements of $M_1$ obtaining the $s$-rigidity of
$\mathcal{M}_1\sqcup\mathcal{M}_2$, too. Thus, the finite value
${\rm max}\{|M_1|+{\rm deg}^{\forall\mbox{\rm -}s}_{\rm
rig}(\mathcal{M}_2),|M_2|+{\rm deg}^{\forall\mbox{\rm -}s}_{\rm
rig}(\mathcal{M}_1)\}$ equals ${\rm deg}^{\forall\mbox{\rm
-}s}_{\rm rig}(\mathcal{M}_1\sqcup\mathcal{M}_2)$.
\endproof

Theorem \ref{th_dis} and Corollary \ref{cor_sem_synt1} immediately
imply:

\begin{corollary}\label{cor_sem_synt1_dis}
For any structures $\mathcal{M}_1$ and $\mathcal{M}_2$ in a
language $\Sigma_1$ of unary predicates the tetrad ${\rm
deg}_4(\mathcal{M}_1\sqcup\mathcal{M}_2)$ has one of the following
possibilities:

$1)$ $(0,0,0,0)$, if $\mathcal{M}_1$ and $\mathcal{M}_2$ are both
semantically and syntactically rigid;

$2)$ $(m,m,n,n)$, if $\mathcal{M}_1$ and $\mathcal{M}_2$ are
finite with $|M_1\,\dot{\cup}\,M_2|=n+1$ elements and some
$\mathcal{M}_i$ is not semantically rigid that is not
syntactically rigid with some minimal $m_1$-elements set
$A_1\subset M_1$ producing ${\rm dcl}(A_1)=M_1$ and some minimal
$m_2$-elements set $A_2\subset M_2$ producing ${\rm
dcl}(A_2)=M_2$, where $m=m_1+m_2\leq n-1$;

$3)$ $(0,\nu,0,\infty)$, if $\mathcal{M}_1\sqcup\mathcal{M}_2$ is
infinite, $\mathcal{M}_1$ and $\mathcal{M}_2$ are semantically
rigid but some of them is not syntactically rigid, with
$1\leq\nu\leq\infty$, $\nu={\rm deg}^{\exists\mbox{\rm -}{\rm
synt}}_{\rm rig}(\mathcal{M}_1)+{\rm deg}^{\exists\mbox{\rm -}{\rm
synt}}_{\rm rig}(\mathcal{M}_2)$;

$4)$ $(\mu,\nu,\infty,\infty)$, if
$\mathcal{M}_1\sqcup\mathcal{M}_2$ is infinite, $\mathcal{M}_1$ or
$\mathcal{M}_2$ is not semantically rigid, $\mathcal{M}_1$ or
$\mathcal{M}_2$ is not syntactically rigid, with
$1\leq\mu\leq\nu\leq\infty$, $\mu={\rm deg}^{\exists\mbox{\rm
-}{\rm sem}}_{\rm rig}(\mathcal{M}_1)+{\rm deg}^{\exists\mbox{\rm
-}{\rm sem}}_{\rm rig}(\mathcal{M}_2)$, $\nu={\rm
deg}^{\exists\mbox{\rm -}{\rm synt}}_{\rm rig}(\mathcal{M}_1)+{\rm
deg}^{\exists\mbox{\rm -}{\rm synt}}_{\rm rig}(\mathcal{M}_2)$.
\end{corollary}

\begin{theorem}\label{th_dis_ind}
For any disjoint predicate structures $\mathcal{ M}_1$ and
$\mathcal{ M}_2$, and a set $A\subseteq M_1\cup M_2$, $${\rm
ind}_{\rm rig}((\mathcal{M}_1\sqcup\mathcal{M}_2)/A)={\rm
max}\{{\rm ind}_{\rm rig}(\mathcal{M}_1/(M_1\cap A)),{\rm
ind}_{\rm rig}(\mathcal{M}_2)/(M_2\cap A)\}.$$
\end{theorem}

Proof. By the definition of disjoint union types in $S^1(A)$ are
locally realized either in $\mathcal{M}_1$ or in $\mathcal{M}_2$.
Moreover, they are forced by their restrictions to $M_1$ or $M_2$.
So algebraic types $p(x)\in S^1(A)$ are defined in $\mathcal{M}_1$
or in $\mathcal{M}_2$ by their restrictions to $M_1\cap A$ and to
$M_2\cap A$. Now we collect possibilities for cardinalities of
sets of realizations of algebraic types in $S^1(M_1\cap A)$ and in
$S^1(M_2\cap A)$. We either choose a maximal natural cardinality
obtaining natural $n={\rm ind}_{\rm
rig}((\mathcal{M}_1\sqcup\mathcal{M}_2)/A)$ with $n={\rm
max}\{{\rm ind}_{\rm rig}(\mathcal{M}_1/(M_1\cap A)),{\rm
ind}_{\rm rig}(\mathcal{M}_2)/(M_2\cap A)\}$ or there are no
maximal natural cardinality with both ${\rm ind}_{\rm
rig}((\mathcal{M}_1\sqcup\mathcal{M}_2)/A)=\omega$ and ${\rm
max}\{{\rm ind}_{\rm rig}(\mathcal{M}_1/(M_1\cap A)),{\rm
ind}_{\rm rig}(\mathcal{M}_2)/(M_2\cap A)\}=\omega$. \endproof

\section{Variations of rigidity for compositions of structures}

Recall the notions of composition for structures and theories.

\medskip
{\bf Definition} \cite{EKS_CT}.  Let $\mathcal{M}$ and
$\mathcal{N}$ be structures of relational languages
$\Sigma_\mathcal{M}$ and $\Sigma_\mathcal{N}$ respectively. We
define the {\em composition} $\mathcal{M}[\mathcal{N}]$ of
$\mathcal{M}$ and $\mathcal{N}$ satisfying the following
conditions:

1)
$\Sigma_{\mathcal{M}[\mathcal{N}]}=\Sigma_\mathcal{M}\cup\Sigma_\mathcal{N}$;

2) $M[N]=M\times N$, where $M[N]$, $M$, $N$ are universes of
$\mathcal{M}[\mathcal{N}]$, $\mathcal{M}$, and $\mathcal{N}$
respectively;

3) if $R\in\Sigma_\mathcal{M}\setminus\Sigma_\mathcal{N}$,
$\mu(R)=n$, then $((a_1,b_1),\ldots,(a_n,b_n))\in
R_{\mathcal{M}[\mathcal{N}]}$ if and only if $(a_1,\ldots,a_n)\in
R_{\mathcal{M}}$;

4) if $R\in\Sigma_\mathcal{N}\setminus\Sigma_\mathcal{M}$,
$\mu(R)=n$, then $((a_1,b_1),\ldots,(a_n,b_n))\in
R_{\mathcal{M}[\mathcal{N}]}$ if and only if $a_1=\ldots =a_n$ and
$(b_1,\ldots,b_n)\in R_{\mathcal{N}}$;

5) if $R\in\Sigma_\mathcal{M}\cap\Sigma_\mathcal{N}$, $\mu(R)=n$,
then $((a_1,b_1),\ldots,(a_n,b_n))\in
R_{\mathcal{M}[\mathcal{N}]}$ if and only if $(a_1,\ldots,a_n)\in
R_{\mathcal{M}}$, or $a_1=\ldots =a_n$ and $(b_1,\ldots,b_n)\in
R_{\mathcal{N}}$.

The theory $T={\rm Th}(\mathcal{M}[\mathcal{N}])$ is called the
{\em composition} $T_1[T_2]$ of the theories $T_1={\rm
Th}(\mathcal{M})$ and $T_2={\rm Th}(\mathcal{N})$.

\medskip
By the definition, the composition $\mathcal{M}[\mathcal{N}]$ is
obtained replacing each element of $\mathcal{M}$ by a copy of
$\mathcal{N}$.

\medskip
{\bf Definition} \cite{EKS_CT}. The composition
$\mathcal{M}[\mathcal{N}]$ is called {\em $E$-definable} if
$\mathcal{M}[\mathcal{N}]$ has an $\emptyset$-definable
equivalence relation $E$ whose $E$-classes are universes of the
copies of $\mathcal{N}$ forming $\mathcal{M}[\mathcal{N}]$.

\begin{note}\label{rem_def_comp} {\rm It is shown in \cite{EKS_CT} that $E$-definable
compositions $\mathcal{M}[\mathcal{N}]$ uniquely define theories
${\rm Th}(\mathcal{M}[\mathcal{N}])$ by theories ${\rm
Th}(\mathcal{M})$ and ${\rm Th}(\mathcal{N})$ and types of
elements in copies of $\mathcal{N}$ are defined by types in these
copies and types for connections between these copies.}
\end{note}

\begin{proposition}\label{prop_wr_comp}
For $E$-definable compositions $\mathcal{M}[\mathcal{N}]$ the
automorphism group ${\rm Aut}(\mathcal{M}[\mathcal{N}])$ is
isomorphic to the wreath product of ${\rm Aut}(\mathcal{M})$ and
${\rm Aut}(\mathcal{N})$:
$$
{\rm Aut}(\mathcal{M}[\mathcal{N}])\simeq{\rm
Aut}(\mathcal{M})\wr{\rm Aut}(\mathcal{N}).
$$
\end{proposition}

Proof. Since all copies of $\mathcal{N}$ are isomorphic in
$\mathcal{M}[\mathcal{N}]$ and form definable $E$-classes each
automorphism $f\in{\rm Aut}(\mathcal{M}[\mathcal{N}])$ is defined
both by the action on the set of $E$-classes, which corresponds to
an automorphism $g\in{\rm Aut}(\mathcal{M})$, and by the the
actions on the $E$-classes, which corresponds to an automorphism
$h$ for copies of $\mathcal{N}$. Therefore $f$ is situated in the
one-to-one correspondence with the pair $(g,h)$ producing a
correspondent element of ${\rm Aut}(\mathcal{M})\wr{\rm
Aut}(\mathcal{N})$. \endproof

\medskip
In view of Remark \ref{rem_def_comp} and Proposition
\ref{prop_wr_comp} we have the following:

\begin{theorem}\label{th_comp1}
For any $E$-definable composition $\mathcal{M}[\mathcal{N}]$ the
following conditions hold: $${\rm deg}^{\exists\mbox{\rm -}{\rm
sem}}_{\rm rig}(\mathcal{M}[\mathcal{N}])={\rm
deg}^{\exists\mbox{\rm -}{\rm sem}}_{\rm rig}(\mathcal{M}),$$ if
$\mathcal{N}$ is semantically rigid, and $${\rm
deg}^{\exists\mbox{\rm -}{\rm sem}}_{\rm
rig}(\mathcal{M}[\mathcal{N}])=|M|\cdot{\rm deg}^{\exists\mbox{\rm
-}{\rm sem}}_{\rm rig}(\mathcal{N}),$$ if $\mathcal{N}$ is not
semantically rigid. In particular, ${\rm deg}^{\exists\mbox{\rm
-}{\rm sem}}_{\rm rig}(\mathcal{M}[\mathcal{N}])$ is finite iff
${\rm deg}^{\exists\mbox{\rm -}{\rm sem}}_{\rm rig}(\mathcal{M})$
and $\mathcal{N}$ are finite, if $\mathcal{N}$ is semantically
rigid, and ${\rm deg}^{\exists\mbox{\rm -}{\rm sem}}_{\rm
rig}(\mathcal{N})$ and $\mathcal{M}$ are finite, if $\mathcal{N}$
is not semantically rigid.
\end{theorem}

Proof. If $\mathcal{N}$ is semantically rigid then it suffices to
find possibilities for automorphisms of $\mathcal{M}$ since in
such a case the semantical rigidity of an inessential expansion of
$\mathcal{M}$ implies the semantical rigidity of correspondent
inessential expansion of $\mathcal{M}[\mathcal{N}]$. Thus, here
${\rm deg}^{\exists\mbox{\rm -}{\rm sem}}_{\rm
rig}(\mathcal{M}[\mathcal{N}])={\rm deg}^{\exists\mbox{\rm -}{\rm
sem}}_{\rm rig}(\mathcal{M})$. If $\mathcal{N}$ is not
semantically rigid then copies of $\mathcal{N}$ in
$\mathcal{M}[\mathcal{N}]$ are automorphically independent, i.e.,
fixing automorphisms for $\mathcal{M}[\mathcal{N}]$ one have to
fix all automorphisms for these copies. Since the smallest set
fixing automorphisms for $\mathcal{N}$ contains ${\rm
deg}^{\exists\mbox{\rm -}{\rm sem}}_{\rm rig}(\mathcal{N})$, we
have at least and minimally at most $|M|\cdot{\rm
deg}^{\exists\mbox{\rm -}{\rm sem}}_{\rm rig}(\mathcal{N})$
elements to fix automorphisms for $\mathcal{M}[\mathcal{N}]$
implying ${\rm deg}^{\exists\mbox{\rm -}{\rm sem}}_{\rm
rig}(\mathcal{M}[\mathcal{N}])=|M|\cdot{\rm deg}^{\exists\mbox{\rm
-}{\rm sem}}_{\rm rig}(\mathcal{N})$. \endproof

\begin{theorem}\label{th_comp2} For any $E$-definable composition $\mathcal{M}[\mathcal{N}]$ the
following conditions hold: $${\rm deg}^{\exists\mbox{\rm -}{\rm
synt}}_{\rm rig}(\mathcal{M}[\mathcal{N}])={\rm
deg}^{\exists\mbox{\rm -}{\rm synt}}_{\rm rig}(\mathcal{M}),$$ if
$N={\rm dcl}(\emptyset)$, and $${\rm deg}^{\exists\mbox{\rm -}{\rm
synt}}_{\rm rig}(\mathcal{M}[\mathcal{N}])=|M|\cdot{\rm
deg}^{\exists\mbox{\rm -}{\rm synt}}_{\rm rig}(\mathcal{N}),$$ if
$N\ne{\rm dcl}(\emptyset)$. In particular, ${\rm
deg}^{\exists\mbox{\rm -}{\rm synt}}_{\rm
rig}(\mathcal{M}[\mathcal{N}])$ is finite iff ${\rm
deg}^{\exists\mbox{\rm -}{\rm synt}}_{\rm rig}(\mathcal{M})$ and
$\mathcal{N}$ are finite, for $N={\rm dcl}(\emptyset)$, and ${\rm
deg}^{\exists\mbox{\rm -}{\rm synt}}_{\rm rig}(\mathcal{N})$ and
$\mathcal{M}$ are finite, for $N\ne{\rm dcl}(\emptyset)$.
\end{theorem}

Proof repeats the proof of Theorem \ref{th_comp1} replacing
automorphism groups by definable closures. \endproof

\medskip
Proposition \ref{prop_rig_sem_synt}, (1), (2) and Theorems
\ref{th_comp1}, \ref{th_comp2} immediately imply:

\begin{corollary}\label{cor_zero}
For any $E$-definable composition $\mathcal{M}[\mathcal{N}]$ and
$s\in\{{\rm sem},{\rm synt}\}$ the following conditions are
equivalent:

$(1)$ ${\rm deg}^{\forall\mbox{\rm -}s}_{\rm
rig}(\mathcal{M}[\mathcal{N}])=0$;

$(2)$ ${\rm deg}^{\forall\mbox{\rm -}s}_{\rm rig}(\mathcal{M})=0$
and ${\rm deg}^{\forall\mbox{\rm -}s}_{\rm rig}(\mathcal{N})=0$.
\end{corollary}

\begin{theorem}\label{th_comp3} For any $s\in\{{\rm sem},{\rm synt}\}$ and $E$-definable composition $\mathcal{M}[\mathcal{N}]$
with $${\rm deg}^{\forall\mbox{\rm -}s}_{\rm
rig}(\mathcal{M}[\mathcal{N}])>0$$ the following conditions are
equivalent:

$(1)$ ${\rm deg}^{\forall\mbox{\rm -}s}_{\rm
rig}(\mathcal{M}[\mathcal{N}])$ is finite;

$(2)$ one of the following conditions hold:

{\rm i)} $\mathcal{M}$ and $\mathcal{N}$ are finite, i.e.
$\mathcal{M}[\mathcal{N}]$ is finite;

{\rm ii)} $\mathcal{M}$ is infinite with ${\rm
deg}^{\forall\mbox{\rm -}s}_{\rm rig}(\mathcal{M})=1$ and ${\rm
deg}^{\forall\mbox{\rm -}s}_{\rm rig}(\mathcal{N})=0$;

{\rm iii)} $\mathcal{M}$ is infinite and $\mathcal{N}$ is finite
with ${\rm deg}^{\forall\mbox{\rm -}s}_{\rm
rig}(\mathcal{M})\in\omega\setminus\{0,1\}$ and ${\rm
deg}^{\forall\mbox{\rm -}s}_{\rm rig}(\mathcal{N})=0$;

{\rm iv)} $\mathcal{M}$ is a singleton and $\mathcal{N}$ is
infinite with ${\rm deg}^{\forall\mbox{\rm -}s}_{\rm
rig}(\mathcal{N})\in\omega\setminus\{0\}$.

Here there are the following possibilities:

{\rm a)} ${\rm deg}^{\forall\mbox{\rm -}s}_{\rm
rig}(\mathcal{M}[\mathcal{N}])=({\rm deg}^{\forall\mbox{\rm
-}s}_{\rm rig}(\mathcal{M})-1)\cdot|N|+1$, if the case {\rm i)} or
{\rm iii)} is satisfied with ${\rm deg}^{\forall\mbox{\rm
-}s}_{\rm rig}(\mathcal{N})=0$;

{\rm b)} ${\rm deg}^{\forall\mbox{\rm -}s}_{\rm
rig}(\mathcal{M}[\mathcal{N}])=(|M|-1)\cdot|N|+{\rm
deg}^{\forall\mbox{\rm -}s}_{\rm rig}(\mathcal{N}),$ if the case
{\rm i)} is satisfied with ${\rm deg}^{\forall\mbox{\rm -}s}_{\rm
rig}(\mathcal{N})>0$;

{\rm c)} ${\rm deg}^{\forall\mbox{\rm -}s}_{\rm
rig}(\mathcal{M}[\mathcal{N}])=1,$ if the case {\rm ii)} is
satisfied;

{\rm d)} ${\rm deg}^{\forall\mbox{\rm -}s}_{\rm
rig}(\mathcal{M}[\mathcal{N}])={\rm deg}^{\forall\mbox{\rm
-}s}_{\rm rig}(\mathcal{N}),$ if the case {\rm iv)} is satisfied.
\end{theorem}

Proof. At first we notice that ${\rm deg}^{\forall\mbox{\rm
-}s}_{\rm rig}(\mathcal{M})>0$ or ${\rm deg}^{\forall\mbox{\rm
-}s}_{\rm rig}(\mathcal{N})>0$ in view of Corollary
\ref{cor_zero}.

Now by the definition $\mathcal{M}[\mathcal{N}]$ is finite iff
$\mathcal{M}$ and $\mathcal{N}$ are finite. In such a case we have
the following possibilities:

$\bullet$  ${\rm deg}^{\forall\mbox{\rm -}s}_{\rm
rig}(\mathcal{M}[\mathcal{N}])=({\rm deg}^{\forall\mbox{\rm
-}s}_{\rm rig}(\mathcal{M})-1)\cdot|N|+1$, if ${\rm
deg}^{\forall\mbox{\rm -}s}_{\rm rig}(\mathcal{N})=0$, since the
rigidity of $\mathcal{M}[\mathcal{N}]$ can be achieved here taking
all elements in ${\rm deg}^{\forall\mbox{\rm -}s}_{\rm
rig}(\mathcal{M})-1$ copies of $\mathcal{N}$ with one additional
element witnessing the degree ${\rm deg}^{\forall\mbox{\rm
-}s}_{\rm rig}(\mathcal{M})$ defining rigidly all $E$-classes for
copies of $\mathcal{N}$ which are rigid by ${\rm
deg}^{\forall\mbox{\rm -}s}_{\rm rig}(\mathcal{N})=0$; it
corresponds the case i) with a);

$\bullet$ ${\rm deg}^{\forall\mbox{\rm -}s}_{\rm
rig}(\mathcal{M}[\mathcal{N}])=(|M|-1)\cdot|N|+{\rm
deg}^{\forall\mbox{\rm -}s}_{\rm rig}(\mathcal{N}),$ if  ${\rm
deg}^{\forall\mbox{\rm -}s}_{\rm rig}(\mathcal{N})>0$, since the
rigidity of $\mathcal{M}[\mathcal{N}]$ can be achieved here taking
all elements in $(|M|-1)$ copies of $\mathcal{N}$ with ${\rm
deg}^{\forall\mbox{\rm -}s}_{\rm rig}(\mathcal{N})$ additional
elements in the last copy of $\mathcal{N}$; it corresponds the
case i) with b).

$(1)\Rightarrow(2)$. Let ${\rm deg}^{\forall\mbox{\rm -}s}_{\rm
rig}(\mathcal{M}[\mathcal{N}])>0$ is finite. We can assume that
$\mathcal{M}$ is infinite or $\mathcal{N}$ is infinite. We have
the following possibilities:

$\bullet$ ${\rm deg}^{\forall\mbox{\rm -}s}_{\rm
rig}(\mathcal{M})=1$ and ${\rm deg}^{\forall\mbox{\rm -}s}_{\rm
rig}(\mathcal{N})=0$, that is any element of
$\mathcal{M}[\mathcal{N}]$ rigidly defines its $E$-class and all
$E$-classes, too, by ${\rm deg}^{\forall\mbox{\rm -}s}_{\rm
rig}(\mathcal{M})=1$, such that all copies of $\mathcal{N}$ in
these $E$-classes are rigid by ${\rm deg}^{\forall\mbox{\rm
-}s}_{\rm rig}(\mathcal{N})=0$; it corresponds the case ii) with
c);

$\bullet$ ${\rm deg}^{\forall\mbox{\rm -}s}_{\rm
rig}(\mathcal{M})\in\omega\setminus\{0,1\}$ and ${\rm
deg}^{\forall\mbox{\rm -}s}_{\rm rig}(\mathcal{N})=0$; here we
require that $\mathcal{N}$ is finite, since otherwise we can take
arbitrary many elements in some $E$-classes which do not imply the
rigidity in view of ${\rm deg}^{\forall\mbox{\rm -}s}_{\rm
rig}(\mathcal{M})\geq 2$; here we have the case iii) with a).

$\bullet$ $\mathcal{M}$ is a singleton and $\mathcal{N}$ is
infinite with ${\rm deg}^{\forall\mbox{\rm -}s}_{\rm
rig}(\mathcal{N})\in\omega\setminus\{0\}$, here ${\rm
deg}^{\forall\mbox{\rm -}s}_{\rm rig}(\mathcal{M})=0$,
$\mathcal{M}[\mathcal{N}]\simeq\mathcal{N}$ and therefore ${\rm
deg}^{\forall\mbox{\rm -}s}_{\rm
rig}(\mathcal{M}[\mathcal{N}])={\rm deg}^{\forall\mbox{\rm
-}s}_{\rm rig}(\mathcal{N})$.

If $\mathcal{N}$ is infinite with ${\rm deg}^{\forall\mbox{\rm
-}s}_{\rm rig}(\mathcal{N})\in\omega\setminus\{0\}$ and
$|\mathcal{M}|\geq 2$ then we can not obtain the rigidity for all
$E$-classes taking arbitrary many elements in some $E$-classes
that contradicts the condition ${\rm deg}^{\forall\mbox{\rm
-}s}_{\rm rig}(\mathcal{M}[\mathcal{N}])\in\omega$.

$(2)\Rightarrow(1)$. Since each finite structure has finite
degrees of rigidity it suffices to show that ${\rm
deg}^{\forall\mbox{\rm -}s}_{\rm rig}(\mathcal{M}[\mathcal{N}])$
is finite if $\mathcal{M}$ is infinite or $\mathcal{N}$ is
infinite with the conditions ii), iii), iv). We observe that ii)
implies c), iii) implies a), and iv) implies d) confirming a
finite value of that degree.
\endproof

\section{Conclusion}

We studied possibilities for the degrees and indexes of rigidity,
both for semantical and syntactical cases. Links of these
characteristics and their possible values are described. We
studied these values and dynamics for structures in some
languages, for some natural operations including disjoint unions
and compositions of structures. A series of examples illustrates
possibilities of these characteristics. It would be interesting to
continue this research describing possible values of degrees and
indexes for natural classes of structures and their theories.

Sergey V. Sudoplatov

Sobolev Institute of Mathematics,

Novosibirsk State Technical
University,

Novosibirsk, Russia

E-mail: sudoplat@math.nsc.ru

\end{document}